\input amstex
\documentstyle{amsppt}
\magnification\magstep1
\nologo
\NoRunningHeads
\NoBlackBoxes
\def\A{\frak A}
\def\Int{\mathop{\roman{Int}}\nolimits}

\def\Fix{\mathop{\roman{Fix}}\nolimits}
\def\conj{\mathop{\roman{conj}}\nolimits}
\def\oo{\varnothing}
\let\ge\geqslant
\let\le\leqslant
\def\C{{\Bbb C}}
\def\R{{\Bbb R}}
\def\Z{{\Bbb Z}}
\def\Cp#1{\C\roman P^{#1}}

\def\Rp#1{\R\roman P^{#1}}
\def\barCP#1{\overline{\C\roman P}^{#1}}
\def\e{\varepsilon}
\def\til#1{\widetilde{#1}}
\let\tm\proclaim
\let\endtm\endproclaim
\def\buildrel#1\over#2{\mathrel{\mathop{\null#2}\limits^{#1}}}
\def\Del{\Delta}
\def\bet{\beta}
\def\alp{\alpha}                
\def\lam{\lambda}               \def\Lam{\Lambda}
\def\gam{\gamma}
\topmatter
\title
On imaginary plane curves and 
Spin quotients of complex surfaces
by complex conjugation
\endtitle
\author
Sergey Finashin and Eugenii Shustin
\endauthor
\rightheadtext{QUOTIENTS OF COMPLEX SURFACES BY COMPLEX CONJUGATION}
\thanks
The second author was partially supported by the grant
no. 6836-1-9 of the Ministry of Science and Technology, Israel
\endthanks
\keywords Imaginary Singularities, Quotients by Complex Conjugation
\endkeywords
\address
Department of Mathematics, Middle East Technical University,
Ankara 06531 Turkey
\endaddress
\email 
serge\,\@\,rorqual.cc.met.edu.tr
\endemail
\address
School of Mathematical Sciences, Tel Aviv University,
Ramat Aviv, Tel Aviv 69978, Israel
\email shustin\,\@\,math.tau.ac.il
\endemail
\endaddress
\subjclass Primary: 14P25, 57N13
\endsubjclass
\abstract
It is proven that for any topological or analytical types of isolated
singular points of plane curves, there exists a non-real 
irreducible plane algebraic curve
of degree $d$ which goes through $d^2$ real distinct points and has 
imaginary singular points of the given types.
This result is used to construct a series of examples of complex algebraic
surfaces $X$ defined over $\R$ whose quotients $Y=X/\conj$ by the complex
conjugation $\conj$ are $Spin$ simply connected 4-manifolds with
signature $16k$, for arbitrary integer $k>0$.
In the previously known examples the signature of
$Spin$ and simply connected quotients $Y$ was zero.
\endabstract 
\endtopmatter

\document
\heading
\S1. Introduction
\endheading

Given a complex algebraic surface $X$ defined over $\R$,
we denote by $X_\R$ the
fixed point set of the complex conjugation $\conj\: X\to X$, and by 
$Y=X/\conj$ the quotient space.
For nonsingular $X$ its quotient
$Y$ is a closed $4$-manifold, which inherits from $X$
an orientation and smooth structure making
the quotient map, $q\:X\to Y$, an 
orientation preserving smooth double covering branched along $X_\R$.
The natural question is to describe
the differential topology types of $4$-manifolds which can arise as
the quotients $Y$ (cf. \cite{D}).

In a plenty of examples the condition that $Y$ is simply connected
(which is the case for example if $X$ is simply connected and $X_\R\ne\oo$)
implies that $Y$ splits into a connected sum of copies of
$\Cp2$, $\barCP2$ and $S^2\times S^2$.
This property for $K3$ surfaces was noticed by S. Donaldson \cite{D}.
S. Akbulut \cite{Ak} set up it for a family
of double planes branched along real curves.
In  \cite{F1} the family of such double planes
was essentially extended and in \cite{F2}  this decomposability property
was set up for
the complete intersection surfaces 
which are constructed by a small perturbation method
(this is a straightforward method which produces however
a plenty of examples of arbitrary multi-degree).
Surprisingly, no examples of complex algebraic
surfaces were known which had
Spin and simply connected quotients $Y$ with non-vanishing signature.
On the other hand,
in certain cases (for example for real double planes $X\to\Cp2$) it is 
not difficult to show that {\it if the quotient $Y$ is
Spin and simply connected then its signature vanish}.

Recently R. Gompf informed one of the authors about the construction
of similar examples in symplectic category.
Namely, by taking fiber sums of elliptic surfaces he constructed
a symplectic $4$-manifold $X$ 
with an anti-symplectic involution $\conj\: X\to X$, which has
simply connected and Spin quotient $Y=X/\conj$ with
signature $\sigma(Y)<0$.
The method of Gompf, probably, can be developed to produce
examples of elliptic surfaces $X$ with the same properties of the 
quotients.

In the present paper we construct examples of a different nature,
when $X$ are  {\it algebraic surfaces of general type\/}
and their quotients $Y$ have {\it positive\/} signature.

\tm{Theorem 1.1} 
There exist  nonsingular complex algebraic surfaces $X$
of general type defined over $\R$ which have simply connected Spin
quotients $Y=X/\conj$ with signature $\sigma(Y)=16k$ for arbitrary 
integer $k\ge1$.
\endtm

These surfaces appear as double coverings over
$\Cp2$ blown up at $k$ 
 pairs of complex conjugated imaginary points, 
$z_i,\bar z_i$, $i=1,\dots,k$.
The branching locus of such a covering is the proper image of a plane
real curve, $A\subset\Cp2$, of even degree, $2d$,
which has non-degenerate singularities of multiplicity $6$
 at the above points $z_i,\bar z_i$ 
and the real part $A_\R=A\cap\Rp2$ consisting
of $d^2$ ovals lying separately on $\Rp2$. 
The curve $A$ is obtained by a perturbation of the curve
$B\cup \conj(B)$, where $B$ is a plane {\it imaginary\/}
(that is not invariant under the complex conjugation)
curve of degree $d$ which
intersects $\Rp2$ precisely at $d^2$ points and has non-degenerate
singularities of multiplicity $6$ at $z_i$, $i=1,\dots,k$
as its only singularities.

Existence of such a curve $B$ follows from the following general result.

\proclaim{Theorem 1.2}
For any given topological or analytic types of isolated singular points
of plane algebraic curves, there exists an imaginary
algebraic curve
in $\Cp2$ of an arbitrarily large degree $d$
which has prescribed numbers of
imaginary singular points of the given types
as its only singularities
and exactly $d^2$ real points.
\endproclaim

We prove Theorem 1.2 in \S2 and Theorem 1.1 in \S3.

\vskip0.5truecm

We thank Prof. A. Libgober (Chicago) for consulting us on the topology
of singular curves on algebraic surfaces.



\head \S 2. Singular plane imaginary curves with maximal number of real points
\endhead

For the sake of simplicity we denote plane curves and 
polynomials defining them by the same symbol.

In this section we prove the following

\proclaim{Theorem 2.1}
For any $k\ge 1$, $n\ge 0$ there exists a curve
in $\Cp2$ of degree
$d=3\cdot 2^{3k+n}$ which has $2^{2n}$
imaginary non-degenerate singular points of
multiplicity $m=2^{k+1}$ and exactly $d^2$ distinct real points.
\endproclaim

This statement immediately implies that
the curve in Theorem 2.1 is irreducible and imaginary,
and at any its real point this curve is non-singular and transversal to
$\Rp2$.
This is because a curve $B\subset\Cp2$
cannot have more then $(\deg B)^2$ real points unless it has
a real irreducible component, which follows from Bezout's theorem
applied to $B$ and $\bar B$. On the other hand, for each real component
$F$ of $B$ the intersection $F\cap\Rp2$ contains either infinite 
number of points, or only isolated (hence singular) points,
which number is bounded from above by
$p=\frac12(\deg F-1)(\deg F-2)$ (arithmetic genus of the curve).

Let us first derive  Theorem 1.2 from
Theorem 2.1.

\demo{Proof of Theorem 1.2}
As it is mentioned above, the real points of the curve $C$ from Theorem 2.1
are transversal to $\Rp2$ and therefore
persist under small deformations. On the other hand,
any topological or analytic singularity is
determined by the $(\mu+1)$-jet of the local Taylor series
at a singular point, where $\mu$ is Milnor number.
It is known that under the condition (see, for instance [DG])
$$(m_1+1)+...+(m_r+1)\le d+1,\tag 2.2$$
for polynomials of degree $d$, the $m_i$-jets, $i=1,...,r$,
of their local Taylor series at any $r$ distinct points can be
prescribed independently. Hence, given certain singularities,
one takes the curve $C$ from Theorem 2.1 with sufficiently
big $k,n$, and then
by a small variation of $C$ in the space of curves
of degree $d$, prescribes the jets,
defining these singularities at some non-degenerate singular points
of $C$. Assume that the curve $\widetilde C$, obtained in the latter
step has extra singular points $w_1,...,w_p$. 
Using independence of jets mentioned above, one deforms $C$ into 
a curve $\widetilde C_i$, which has the same jets at the points
giving prescribed singularities and which does not pass through $w_i$,
($1\le i\le p$). Then by Bertini's theorem, one kills all extra
singularities by variation in the linear system
$$\lam\widetilde C+\lam_1\widetilde C_1+...+\lam_p\widetilde C_p,\quad
\lam,\lam_1,...,\lam_p\in\R\quad\qed$$
\enddemo

\demo{Proof of Theorem 2.1} We proceed by induction on $k$ in the case
$n=0$, and then by induction on $n$ with a fixed $k\ge1$.

{\it Step 1}. Let $n=0$, $k=1$, $d=3$, $m=2$. Let $C_3$ 
be a real cuspidal
cubic curve and $C'_3$ be a cubic curve, meeting $C_3$ at
9 distinct real points (for instance, we can take $C'_3$ to be the
union of three real straight lines intersecting $C_3$ each at three
real points). The pencil
$$\Lam=\{\alp C_3+\bet C'_3\ |\ [\alp:\bet]\in\Rp1\}$$
intersects the discriminant in the space $\Rp9$ of real cubics
at $C_3$ with multiplicity 2, because cuspidal curves
form a cuspidal edge of the discriminant (see, for instance [Gu, \S 2]).
Since a small variation of the pencil $\Lam$ is a pencil with 9 real
base points as well, we can move $\Lam$ in such a position that
its complexification will meet the discriminant at two imaginary
points near $C_3$, which correspond to imaginary nodal cubics
with 9 real points.

{\it Step 2}. The induction procedure in the case $n=0$ will be
as follows:
given an imaginary
curve of degree $d$ with a non-degenerate
point of multiplicity $m<d$ and $d^2$ real points,
we construct a curve of degree $8d$ with a non-degenerate
point of multiplicity $2m$ and $64d^2$ real points.

Let $C_d$ be an imaginary
curve of degree $d$ with
an imaginary
non-degenerate singular point $z$ of multiplicity $m$, having exactly
$d^2$ real points, $w_i$, $i=1,...,d^2$.
It is easily seen, say, by (2.2), that for polynomials of degree $d$,
the $m$-jet 
at $z$ and $0$-jet at any point $z'\ne z$
are independent. Hence, varying a little
the curve $C_d$ we can provide that none of the points 
$w_i$, $1\le i\le d^2$, lies on the
real straight line $L$ through $z$, and $L$ is not tangent to $C_d$ at $z$.
Let us choose real
projective coordinates $[x_0\!:\!x_1\!:\!x_2]$ so that $L=\{x_0=0\}$, and
$$w_i=[1\!:\!x_i\!:\!y_i],\quad x_i>0,\ y_i>0,\quad i=1,...,d^2.$$
The transformation $[x_0\!:\!x_1\!:\!x_2]\mapsto[x_0^2\!:\!x_1^2\!:\!x_2^2]$ turns
the points $w_1,...,w_{d^2}$ into $4d^2$ real points of the curve
$$C_{2d}(x_0,x_1,x_2)\buildrel{\text{def}}\over{=}C_d(x_0^2,x_1^2,x_2^2)
\tag 2.3$$
of degree $2d$. The point $z$ will turn 
into two imaginary singular points
$z_1,z_2$ of $C_{2d}$, such that $z_i,i=1,2,$ is a center of $m$ 
non-singular
branches quadratically 
tangent to the straight line $L_i$ through
$z_i$ and $[1\!:\!0\!:\!0]$: indeed, in affine coordinates $x=x_0/x_2,y=x_1/x_2$
a local branch of $C_d$ centered at $z=(0,a)$, $a\ne 0$, has equation
$$y=a+bx+O(x^2),\quad b\ne 0,$$
that is transformed into
$$y^2=a+bx^2+O(x^4),$$
which defines two non-singular branches
$$y=\pm\sqrt{a}\left(1+{b\over 2a}x^2\right)+O(x^3),$$
centered at the points $z_{1,2}=(0,\pm\sqrt{a})$ and
quadratically tangent to the lines $y=\pm\sqrt{a}$, 
respectively.

{\it Step 3}. Real conics tangent to $L_1$ at $z_1$, form a pencil. 
Let
$K$ be a generic conic in this pencil, having a real oval.
We choose on $K$ three generic real points and change coordinates in such
a way that these points become 
$[1\!:\!0\!:\!0],[0\!:\!1\!:\!0],[0\!:\!0\!:\!1]$,
respectively. The Cremona transformation (see detail in [Wa])
$$[x_0\!:\!x_1\!:\!x_2]\mapsto[x_1x_2\!:\!x_0x_2\!:\!x_0x_1]\tag 2.4$$
takes the conic $K$ to a real straight line $L^*$, and the curve
$C_{2d}$ to a curve $C_{4d}$ of degree $4d$ having
$4d^2$ real non-singular points $v_i,i=1,...,4d^2$,
three non-degenerate singularities of multiplicity $2d$ at the real points
$[1\!:\!0\!:\!0],[0\!:\!1\!:\!0],[0\!:\!0\!:\!1]$, and an imaginary
singular point 
$z_1^*\in L^*$ (coming from $z_1$),
which is a center of $m$ non-singular local branches
quadratically tangent to $L^*$.

{\it Step 4}.
Keeping singularity at $z_1^*$, we shall deform
the ordinary singular points $[1\!:\!0\!:\!0],[0\!:\!1\!:\!0],[0\!:\!0\!:\!1]$ of $C_{4d}$ in
order to get $12d^2$ more real non-singular points. We apply
a modified Viro method (see the original Viro method in [Vi1, Vi2, Vi3, Ri]
and [GKZ, Chapter 11],
and relevant modifications in [Sh]).

In the affine coordinates $x=x_1/x_0$,
$y=x_2/x_0$ the curve $C_{4d}$ is defined by a polynomial
$F_0(x,y)$ with Newton triangle
$$\Del_0=\text{conv}((2d,0),\ (0,2d),\ (2d,2d)).$$
Consider the polynomials
$$\eqalign{&F_1(x,y)=x^{2d}y^{2d}F_0(y^{-1},x^{-1}),\cr
&F_2(x,y)=x^{-2d}y^{4d}F_0(x,xy^{-1}),\cr
&F_3(x,y)=x^{4d}y^{-2d}F_0(x^{-1}y,y).\cr}$$
It is easily seen that
\roster
\item"(i)" these polynomials have Newton triangles
$$\eqalign{&\Del_1=\text{conv}((0,0)\ (2d,0),\ (0,2d)),\cr
&\Del_2=\text{conv}((0,2d),\ (0,4d),\ (2d,2d)),\cr
&\Del_3=\text{conv}((2d,0),\ (4d,0),\ (2d,2d)),\cr}$$
respectively, which together with $\Del_0$ form a subdivision
of the Newton triangle $T_{4d}=$conv$((0,0),(4d,0),(0,4d))$ of a generic
polynomial of degree $4d$,
\item"(ii)" each curve $F_i$ has exactly $4d^2$ real points in
$(\R^*)^2$, $i=0,1,2,3$,
\item"(iii)" the coefficients of any common monomial in $F_i,F_j$, $i\ne j$,
coincide.
\endroster
In particular, the polynomials $F_0,F_1,F_2,F_3$ uniquely define a set
of coefficients $a_{ij}$ of monomials $x^iy^j$, $i+j\le 4d$.
Let $\nu\:T_{4d}\to\R$ be a convex piecewise-linear function, whose linearity
domains are just $\Del_0,\Del_1,\Del_2,\Del_3$, which vanishes on
$\Del_0$ and is integral-valued at integral points.
We will look for the desired deformation in the Viro-type family
$$\Phi_t(x,y)=\sum_{i+j\le 4d}A_{ij}(t)x^iy^jt^{\nu(i,j)},\quad t>0,
\tag 2.5$$
where $A_{ij}(t)=a_{ij}+O(t)$, $i+j\le 4d$, are smooth functions.
We find out $A_{ij}(t)$ assuming that the curves
$\Phi_t$ have $m$ smooth branches at $z^*_1$ tangent to $L^*$.
This means that in suitable coordinates in a neighborhood of $z^*_1$
coefficients of monomials lying under the Newton diagram
$[(0,m),(2m,0)]$ vanish, that imposes $m(m+1)$ linear conditions
on coefficients of the polynomial $\Phi_t(x,y)$:
$$\sum_{i+j\le 4d}\alp^{(k)}_{ij}A_{ij}(t)t^{\nu(i,j)}=0,\quad
k=1,...,m(m+1).\tag 2.6$$
These $m(m+1)$
linear equations, imposed by singularity at $z^*_1$,
are independent in the space of
polynomials with Newton triangle $\Del_0$, because the Cremona transformation
(2.4) reduces this to the evident independence in the space of
polynomials of degree $2d$. In particular, there exists a set
$$\Lam\subset\Del_0\cap\Z^2,\quad \text{card}(\Lam)=m(m+1),$$
such that
$$\det\left(\alp^{(k)}_{ij}\right)_{k=1,...,m(m+1)\atop(i,j)\in\Lam}\ne 0.$$
Since $\nu\Big|_{\Del_0}=0$ and $A_{ij}=a_{ij},i+j\le 4d$, $t=0$ is
a solution of (2.6), the latter inequality by the implicit function
theorem provides existence 
of a solution $A_{ij}(t)=a_{ij}+O(t)$,
$i+j\le 4d$, of (2.6) for small $t>0$.

{\it Step 5}. Now we shall show that the curve $\Phi_t$ has
$16d^2$ real points. Let $Q\subset(\R^*)^2$ be a compact
neighborhood of all the real points of the curves $F_i$,
$i=0,1,2,3$. By construction
$\Phi_t(x,y)=F_0(x,y)+O(t)$; hence for a sufficiently small $t>0$, 
the curve
$\Phi_t$ has in $Q$ at least $4d^2$ real points close to these
of the curve $F_0$. Let $l_k(i,j)=\gam_{0k}+\gam_{1k}i
+\gam_{2k}j$ be the linear function equal to $\nu(i,j)$ on $\Del_k$,
$1\le k\le 3$. The substitution of $\nu_k=\nu-l_k$ for $\nu$ in
(2.5) gives us the polynomial
$$\Phi_{t,k}(x,y)=\sum_{i+j\le 4d}A_{ij}(t)x^iy^jt^{\nu_k(i,j)}=
F_k(x,y)+O(t);$$
hence for a sufficiently small $t>0$
the curve $\Phi_{t,k}$ has
in $Q$ at least $4d^2$ real points close to these for the
curve $F_k$. On the other hand, substitution 
of $\nu_k$
for $\nu$ is equivalent to multiplication 
of the polynomial
by a positive number and the coordinate change
$$T_k(x,y)=(xt^{-\gam_{1k}},yt^{-\gam_{2k}}).$$
Therefore the curve $\Phi_t$ has at least $4d^2$ real points in
$T_k^{-1}(Q)$, $k=1,2,3$. Observing that 
$Q$, $T_1^{-1}(Q)$, $T_2^{-1}(Q)$, $T_3^{-1}(Q)$ are disjoint for small
$t>0$, and that $16d^2$ is the maximal possible number of real
points for an irreducible imaginary curve of degree $4d$,
one concludes that the curve $\Phi_t$ has exactly
$16d^2$ real points.

{\it Step 6}. Now let us make a real coordinate change in the
projective plane such that $L^*=\{x_0=0\}$ and all $16d^2$
real points of the projective closure $\hat\Phi_t$ of the
curve $\Phi_t$ lie in the
domain $x_0=1,\ x_1>0,\ x_2>0$. Consider the curve
$$C_{8d}(x_0,x_1,x_2)\buildrel{\text{def}}\over{=}\hat
\Phi_t(x_0^2,x_1^2,x_2^2)$$
of degree $8d$. Each real point of the curve $\hat\Phi_t$ turns into
4 real points of the curve $C_{8d}$, and the singular point $z^*_1$ of
$\hat\Phi_t$ 
turns into two ordinary singular points of
multiplicity $2m$: indeed, in the coordinates $x=x_1/x_2,y=x_0/x_2$
the equation
$$y=b(x-a)^2+O((x-a)^3),\quad a,b\ne 0,$$
of a local branch of $\hat\Phi_t$ centered at $z_1=(a,0)$, 
turns into
the equation
$$y^2=b(x^2-a)^2+O((x^2-a)^3),$$
which defines two pairs of non-singular transversal branches
centered at the points $(\sqrt{a},0)$, $(-\sqrt{a},0)$.

Hence $C_{8d}$ has $64d^2=(8d)^2$
real singular points and an
imaginary singular point of multiplicity $2m$. This completes the induction
step, and thereby the proof of Theorem 2.1 for $n=0$.

{\it Step 7}. Now we use induction on $n$, applying the transformation
(2.3), so that each time the coordinate system is chosen in order to
provide $x_1/x_0>0$, $x_2/x_0>0$ for all real points of $C_d$, and
$x_0x_1x_2\ne 0$ for all imaginary non-degenerate singular points of $C_d$.
$\qed$\enddemo




\heading 
\S3. Proof of Theorem 1.1
\endheading

Consider an imaginary curve $B\subset\Cp2$ of degree $d$ which
intersects $\Rp2$ in $d^2$ points, $p_1,\dots,p_{d^2}$,
and has non-degenerated singularities of multiplicity $m=6$
at $z_i,\in\Cp2-\Rp2$, $i=1,\dots k$, as its only singularities,
(such $B$ exists by Theorem 1.2).
Then $A_0=B\cup\conj(B)$ is a degree $2d$ curve, which has in addition to
the singularities at $z_i$ and $\bar z_i=\conj(z_i)$, $i=1,\dots,k$,
nodes at $p_j$, $j=1,\dots,d^2$.
We can perturb $A_0$ preserving its singularities at $z_i,\bar z_i$ 
and making the nodes smooth
so that each node gives rise to an oval.
If $B$ is defined by equation $f=0$, then $A_0$ is defined by
$f\centerdot\bar f=0$, where $\bar f$ is the polynomial conjugated to $f$,
and the perturbed curve, which we denote by $A_\e$, is defined by
the polynomial
$g_\e=f\centerdot\bar f-\e h^2$, where  $\e>0$ is a sufficiently small
real number and $h$ is a degree $d$ real homogeneous polynomial
which defines a curve not containing the points $p_i$, $i=1,\dots,d^2$,
and containing $z_i, \bar z_i$, $i=1,\dots,k$, with multiplicity $>3$
(such $h$ can be easily constructed since $d>>k$).

Consider blow up $P\to \Cp2$ centered at points 
$z_i,\bar z_i$, $i=1,\dots,k$, and
denote by $\til{A}_\e\subset P$ the proper image of $A_\e$
and by $\til{B}$ the proper image of $B$.
$\til A_\e$ is a nonsingular divisible by $2$ divisor in $P$;
therefore there exists an algebraic surface $X$ with a double
covering $p\:X\to P$ branched along $\til A_\e$. 
It is well known \cite{A} that there exist two liftings to $X$ of
the complex conjugation $\conj$, induced on $P$ from $\Cp2$.
These liftings are anti-holomorphic involutions,
which commute and give in product
the covering transform of $p$.
The fixed point set of one
of these involution consists of $d^2$ sphere components, which
are projected by $p$ onto the union $W^+$ of $d^2$ discs bounded in $\Rp2$
by the ovals of $A_\R=A\cap\Rp2$.
We choose and 
denote by $\conj_X$ the other involution,
whose fixed point set, $X_\R=\Fix(\conj_X)$, is projected into
the complementary part, $W^-=\Rp2-\Int(W^+)$. 
Denote by $Q=P/\conj$, $Y=X/\conj_X$
the quotients and by $q\:P\to Q$, $q'\:X\to Y$ the quotient maps. 

Arnold noticed in \cite{A}
(in a different but analogous case of double planes)
that factorization by $\Z/2\times\Z/2$ in different order gives the
following commutative diagram
$$
\CD
X     @>q'>>   Y\\
@VpVV         @VVp'V\\
P   @>q>>     Q  ,
\endCD
$$
where the map $p'$ turns out to be
the double covering branched along the surface
$\A_\e=q(W^+)\cup q(\til A_\e)$ called the Arnold surface.

Recall now the well known diffeomorphism
$\Cp2/\conj\cong S^4$ (see, e.g., \cite{M}) which implies
$Q\cong\#_k\barCP2$. Denote by $E_i',E_i''\subset P$ the exceptional 
curves over $z_i,\bar z_i$, $i=1,\dots,k$, and take as the generators of 
$H_2(Q)$ the fundamental classes $[E_i]$ of
$E_i=q(E_1)=q(E_2)$ with the orientations coming from $E_i$.

Further, note that varying $\e\to0$ we get an isotopy in $Q$ between
$\A_\e$ and $q(\til B)$  (which pinch the ovals of $\A_\e$),
 cf. \cite{F1}. This implies that
$\A_\e$ endowed with the orientation coming from $B$ realizes class
$[\A_\e]=6\sum_{i=1}^k[E_i]\in H_2(Q)$.
Applying Hirzebruch signature formula to the branched covering
$Y\to Q$ we obtain 
$$
\sigma(Y)=2\sigma(Q)-\frac12 \A_\e\circ \A_\e=
-2k+\frac126^2k=16k
$$
The formula for the Stiefel--Whitney class $w_2$ of the branched covering
gives
$$
w_2(Y)={p'}^*(w_2(Q))+{p'}^*(\frac12[\A_\e])_2=
{p'}^*(\sum_{i=1}^k[E_i]_2)+{p'}^*(3\sum_{i=1}^k[E_i]_2)=0,
$$
where subscript ``$2$'' denotes the reduction of an integer class 
modulo $2$.

To finish the proof of properties of $Y$ we need only
to notice that it is simply connected.
To see it note first that
the fundamental group $\pi_1(P-\til A)$ is abelian, which follows
for instance from the Nori's theorem \cite{N, Theorem II}. Thus
it is cyclic with the generator represented by a loop around
$\til A$, which follows from the following diagram
$$
\CD
@. 			H^2(\til A) 	@=\Z\\
@. 			@VVV \\
H_1(P-\til A) @>{\cong}>> H^3(P,\til A) \\
@.			@VVV\\
@.			 H^3(P)		@=0\\
\endCD
$$
Therefore $\pi_1(X-p^{-1}(\til A))$ is also cyclic with the generator
represented by a loop around $p^{-1}(\til A)$. Thus $X$ is simply
connected and, therefore $Y$ is simply connected, since $X_\R\ne\oo$.

Finally, note that involution $\conj_X$ is the involution of 
complex conjugation for some embedding $X\to\Cp{N}$, because
as it was noticed by Comessatti (cf. \cite{CP}),
 {\it any\/ anti-holomorphic involution
on a complex algebraic variety $X$ is the restriction of the
complex conjugation on $\Cp{N}$
for some embedding $X\to\Cp{N}$}
(such an embedding is defined by the line bundle $L\otimes\conj_X^*(L)$
for any very ample line bundle $L$ on $X$).
\qed




\enddocument